\documentclass{amsart}

\theoremstyle{definition}

\theoremstyle{remark}
\newtheorem{remark}{Remark}

\usepackage{multicol}
\usepackage{euscript}
\usepackage{amssymb}
\usepackage{amsmath}
\usepackage{amsthm}
\usepackage[dvips]{graphicx}
\usepackage{psfrag}

\newcommand{\ee}{\varepsilon}

\newcommand{\ue}{u^{\varepsilon}}

\begin{document}
\bibliographystyle{plain}
\title[Viscous profiles and singular ODEs]{A connection between viscous profiles and singular ODEs}

\author{Stefano Bianchini}
\address{S.B.: SISSA, via Beirut 2-4 34014 Trieste, Italy}
\email{bianchin@sissa.it}

\author{Laura V. Spinolo}
\address{L.V.S.: Centro De Giorgi, Collegio Puteano, Scuola Normale Superiore,
Piazza dei Ca\-va\-lie\-ri 3, 56126 Pisa, Italy}
\email{laura.spinolo@sns.it}


\begin{abstract}
We deal with the viscous profiles for a class of mixed hyperbolic-parabolic systems. 
We focus, in particular, on the case of the compressible Navier Stokes equation in one space variable written in Eulerian coordinates.
We describe the link between these profiles and a singular ordinary differential equation in the form
\begin{equation}
\label{e:odesing:intro}
         \frac{dV}{dt} =  \frac{1}{  \zeta (V)}  F(V).
\end{equation}
Here $V \in \mathbb{R}^d$ and the function $F$ takes values 
into $\mathbb{R}^d$ and is smooth. The real valued function $\zeta $ is  as well regular: the equation is singular in the sense that 
$\zeta (V)$ can attain the value $0$. 
\end{abstract}
\maketitle 

We focus on mixed hyperbolic-parabolic systems in one space variable in the form
\begin{equation}
\label{e:pare:vispro}
        E(u) u_t + A(u, \, u_x) u_x = B (u) u_{xx}.
\end{equation}
Here the function $u$ takes values in $\mathbb{R}^N$ and depends on the two scalar variables $t$ and $x$. We focus on the case the matrix $B$ is singular, namely
its rank is strictly smaller than $N$. 
In particular,  a conservative system  
$$
   u_t + f(u)_x = \Big( B(u) u_x \Big)_x 
$$
can be written in the form \eqref{e:pare:vispro}. Indeed, one can set 
$$
    A(u, \, u_x ) = D \, f  (u)- B(u)_x, 
$$
where $D \, f$ denotes the jacobian matrix of $f$. In the following we will consider explicitly the case of the compressible Navier Stokes equation in one space variable:
\begin{equation}
\label{e:ns:eul}
       \left\{
       \begin{array}{lll}
              \rho_t + ( \rho v )_x =0 \\
	      (\rho v)_t + \Big( \rho v^2 + p \Big)_x = \displaystyle{ \Big( \nu v_x  \Big)_x } \\
	      \displaystyle{ \Big( \rho e + \rho \frac{v^2}{2}\Big)_t + \Big(v \Big[ \frac{1}{2} \rho v^2 
	      + \rho e + p \Big] \Big)_x = \Big( k \theta_x + 
	      \nu v v_x \Big)_x}. \\
       \end{array}
       \right.
\end{equation}
Here the unknowns $\rho(t, \, x), \, v(t, \, x)$ and $\theta(t, \, x)$ are the density of the fluid, the velocity of the particles in the fluid and the absolute temperature respectively. The function $p= p(\rho, \, \theta) >0$ is the pressure and satisfies $p_{\rho} >0$, while $e$ is the internal energy. In the following we will focus on the case of a polytropic gas, so that $e$ satisfies
\begin{equation}
\label{e:pol}
   e = \frac{R   \theta }{ \gamma -1 },  
\end{equation}
where $R$ is the universal gas constant and $\gamma >1$ is a constant specific of the gas. Finally, by $\nu(\rho)>0$ and $k(\rho)>0$ we denote the viscosity and the heat conduction coefficients respectively.

In the following, we assume that system \eqref{e:pare:vispro} satisfies a set of hypotheses introduced by Kawashima and Shizuta in 
\cite{KawShi:normal}. These conditions were defined
relying on examples with a physical meaning, in particular, up to a change in the dependent variables, they are satisfied by the equations of the hydrodynamics.   
For completeness, we recall these conditions here. 
First, we assume that the  rank of the matrix $B(u)$ is constant and we denote it  by $r$. Also, $B(u)$  
       admits the block decomposition
       \begin{equation}
       \label{e:blockdec:h}
             B(u) =
             \left(
             \begin{array}{cc}
                   0 & 0 \\
                   0 & b(u) \\
             \end{array}
             \right).
       \end{equation}
       The block $b(u)$ belongs to $ \mathbb{M}^{r \times r}$ and there exists a constant $c_b > 0$  such that for every $\vec{\xi} \in \mathbb{R}^r$ 
       \begin{equation}
       \label{e:bpos:vispro}
              \langle b(u) \vec{\xi}, \, \vec{\xi} \rangle
              \ge c_b |\vec{\xi}|^2.
       \end{equation}
       In the previous expression, $\langle \cdot, \, \cdot \rangle$ denotes the standard scalar product in $\mathbb{R}^r$. 

Also, we assume that for every $u$
       the matrix $A(u, \, \, 0)$ is symmetric. We denote by 
       \begin{equation}
       \label{e:blocka:hyp}
              A(u, \, u_x) =
              \left(
             \begin{array}{cc}
                   A_{11}(u) & A_{21}^T(u) \\
                   A_{21}(u) & A_{22}(u, \, u_x) \\
             \end{array}
             \right) \qquad 
              E(u) =
              \left(
             \begin{array}{cc}
                   E_{11}(u) & E_{21}^T(u) \\
                   E_{21}(u) & E_{22}(u) \\
             \end{array} 
             \right)    
           \end{equation}
      the block decomposition of $A$ and $E$ corresponding to \eqref{e:blockdec:h}. Note that only the block $A_{22}$ can depend on $u_x$. Finally, 
            we assume that   
       for every $u$, the matrix
       $E(u)$
       is  symmetric and positive definite.

In the following, we focus on two classes of solutions of \eqref{e:pare:vispro}: travelling waves and steady solutions. A travelling wave is a one variable function
$U(y)$ satisfying 
\begin{equation}
\label{e:tw:vispro} 
         \big[ A ( U, \, U') - \sigma E (U) \big] U ' = B (U) U '' .
\end{equation}
In the previous expression, $\sigma$ is  a real number and is the speed of the wave. From \eqref{e:tw:vispro} one obtains a solution of \eqref{e:pare:vispro} by setting $u(t, \, x) = U (x - \sigma t)$. Steady solutions are solutions of \eqref{e:pare:vispro} that do not depend on time: namely, they are one variable functions $U(x)$ satisfying 
\begin{equation}
\label{e:steady:vispro}
   A(U, \, U') U' = B(U) U''. 
\end{equation}
We also require that $U$ is bounded on $x \in [0, \, + \infty [$ and admits a limit as ${x \to + \infty}$. Steady solutions in this class are sometimes called {\it boundary layers}. In the applications, it is often interesting to focus on the case the speed $\sigma$ in \eqref{e:tw:vispro} is close to an eigenvalue of $E^{-1}A$. Since in general $0$ is not an eigenvalue of $E^{-1}A$, it is useful to distinguish between~\eqref{e:tw:vispro} and \eqref{e:steady:vispro}. 

Travelling waves and steady solutions are powerful tools to study the parabolic approximation of hyperbolic problems: since the literature concerning this topic is extremely big, we just refer to the books by Dafermos \cite{Daf:book} and by Serre \cite{Serre:book} and to the rich bibliography contained therein. Concerning the analysis of the viscous profiles, we refer to  Benzoni-Gavage, Rousset, Serre and Zumbrun \cite{BenRouSerreZum}, to Zumbrun \cite{Zum:review} and to the references therein. 
Loosely speaking, the problem of the parabolic approximation of an hyperbolic system is the following: consider the family of systems 
\begin{equation}
\label{e:app:vispro}
        E(\ue) \ue_t + A(\ue, \, \ee \ue_x) \ue_x = \ee B (\ue) \ue_{xx},
\end{equation}
which reduces to \eqref{e:pare:vispro} via the change of variable $u(t, \, x)= \ue ( \ee t, \, \ee x)$. Letting $\ee \to 0^+$, equation \eqref{e:app:vispro} formally reduces to 
\begin{equation}
\label{e:hyp:vispro}
    E (u) u_t + A(u, \, 0) u_x =0.
\end{equation}
As we pointed out before, the compressible Navier Stokes equation in one space variable can be written in a form like \eqref{e:app:vispro}: in this case, system \eqref{e:hyp:vispro}, which formally is obtained by setting $\ee=0$, is the Euler equation.   
The proof of the convergence of $\ue$ to a solution of \eqref{e:hyp:vispro} is an open problem in the most general case, but convergence results are available in more specific situations (see e.g. Bianchini and Bressan \cite{BiaBrevv}, Ancona and Bianchini \cite{AnBiapro} and  Gisclon \cite{Gisclon:etudes}). Travelling waves have been exploited to study the parabolic approximation \eqref{e:app:vispro} and its relations with the limit \eqref{e:hyp:vispro}, especially to describe such phenomena as the formation of singularities. The study of steady solutions has provided useful information to the analysis of \eqref{e:app:vispro} especially in the case of initial boundary value problems. In the following we will be mainly concerned with viscous profiles having sufficiently small total variation.  

In a previous work~\cite{BiaSpi:rie}, the authors studied the approximation~\eqref{e:app:vispro} for an initial boundary value problem with special data. Concerning the structure of the matrices $E$, $A$ and $B$, we introduced a new condition of block linear degeneracy, which is the following. Let the block $A_{11}(u)$ and $E_{11} (u) $ be as in \eqref{e:blocka:hyp}. The block linear degeneracy precribes that for every real number $\sigma$ the dimension of the kernel of $[A_{11}(u) -\sigma E_{11}(u)]$ is constant with respect to $u$. In other words, the condition of block linear degeneracy says that such a dimension can in general vary when $\sigma$ varies, but it cannot vary when $u$ varies. 

If the condition of block linear degeneracy is violated, then one might face ``pathological" behaviors, in the following sense.  In \cite[Section 2]{BiaSpi:rie} we discuss an example of a system that satisfies all the Kawashima Shizuta hypotheses, but does not satisfy the block linear degeneracy. We exhibit a steady solution \eqref{e:steady:vispro} of this system 
which is not continuously differentiable. Note that in this case it is not completely clear, a priori, what we mean by solution of \eqref{e:steady:vispro}, because we are dealing with non regular functions. The details are given in \cite{BiaSpi:rie}.

On the other side, imposing the block linear degeneracy is restrictive in view of some applications. In particular, as pointed out by Fr{\'e}d{\'e}ric Rousset in~\cite{Rousset:ns}, this condition is satisfied by the compressible Navier Stokes equation written in Lagrangian coordinates, but is violated by the same  equation written in Eulerian coordinates.

The problem is the following. Consider the Navier Stokes equation written in Lagrangian coordinates: 
\begin{equation}
\label{e:ns:lag}
       \left\{
       \begin{array}{lll}
              \tau_t - v_x =0 \\
	      v_t + p_x = \displaystyle{ \Big( \frac{\nu}{\tau} v_x  \Big)_x } \\
	      \displaystyle{ \Big( e +  \frac{v^2}{2}\Big)_t + \Big( p v \Big)_x = \Big( \frac{k}{\tau} \theta_x + 
	      \frac{\nu}{\tau} v v_x \Big)_x}. \\
       \end{array}
       \right.
\end{equation}
Here the unknowns are $\tau(t, \, x)$, $v(t, \, x)$ and $\theta (t, \, x)$: $\tau$ is the \emph{specific volume}, and it is the inverse of the density, $\tau = 1 / \rho$. As in~\eqref{e:ns:eul}, $v$ and $\theta$ denote the velocity and the absolute temperature of the fluid. As pointed out for example in Rousset \cite{Rousset:char}, the compressible Navier Stokes equation written in Lagrangian coordinates can be written in a form like \eqref{e:pare:vispro}, with $E$, $A$ and $B$ satisfying all the Kawashima Shizuta~\cite{KawShi:normal} conditions. In particular, the rank of the viscosity matrix is constantly equal to two. Let $A_{11}$ be the same block as in \eqref{e:blocka:hyp}. In the case of the compressible Navier Stokes written in Lagrangian coordinates $A_{11}$ is a real valued function and it is actually identically equal to $0$. 
This implies, in particular, that the condition of block linear degeneracy is satisfied by the compressible Navier Stokes written in Lagrangian coordinates. Indeed, let $E_{11}(u)$ as in \eqref{e:blocka:hyp}: in this case, $E_{11}(u)$ is a strictly positive, real valued function. It follows that, for every real number $\sigma$, the dimension of the kernel of $- \sigma E_{11}(u)$ does not depend on $u$. 

By direct computations one can then verify that the viscous profiles of the Navier Stokes equation written in Lagrangian coordinates satisfy an ordinary differential equation which does not have the singularity exhibited by~\eqref{e:odesing:intro}.

Let us consider now~\eqref{e:ns:eul}, the Navier Stokes equation written in Eulerian coordinates. We want to write it in the form 
\begin{equation}
\label{e:ns:abstract}
    E (u) u_t + A (u, \, u_x )u_x = B(u) u_{xx},
\end{equation}
requiring that the matrix $A(u, \, 0)$ is \emph{symmetric}. In the following, we assume that $\rho$ is bounded away from $0$, say $\rho \ge c_{\rho} >0$ for a suitable constant $c_{\rho}$. This implies that the system does not reach the vacuum.  

We proceed as in Kawashima and Shizuta~\cite{KawShi:normal} and by multiplying~\eqref{e:ns:eul} on the left by a suitable nonsingular matrix we eventually obtain that~\eqref{e:ns:eul} can be written in the form \eqref{e:ns:abstract} for 
\begin{equation}
\label{e:h:ns:b}
           E (\rho, \, v, \, \theta) =  \frac{1}{\theta}
           \left(
           \begin{array}{ccc}
                                   p_{\rho} / \rho & 0 & 0 \\
                                   0 & \rho   & 0   \\
                                   0 &  0 & \rho e_{\theta} / \theta \\ 
           \end{array}
           \right)
           \qquad 
          B (\rho, \, v, \, \theta) =  \frac{1}{\theta}
           \left(
           \begin{array}{ccc}
                                 0 & 0 & 0      \\
                                   0 & \nu  &  0  \\
                                   0 &  0 &   k / \theta \\ 
           \end{array}
           \right)           
\end{equation}
and 
\begin{equation}
\label{e:h:ns:a}
       A (\rho, \, v, \, \theta, \, \rho_x, \, v_x, \, \theta_x ) =   
       \frac{1}{\theta}
        \left(
           \begin{array}{ccc}
                     p_{\rho} v / \rho    &  
                    p_{\rho}  & 0 \\
                    p_{\rho} & 
                    \rho v -  \nu' \rho_x &  p_{\theta} \\
                     0 & 
                      p_{\theta}  - 
                       \nu v_x /  \theta  
                      &    \rho v e_{\theta} / \theta  
                       -  k' \rho_x  /  \theta  \\ 
           \end{array}
           \right)
\end{equation}
In the previous expression, we denote by $p_{\rho}$ and $p_{\theta}$ the partial derivative of $p$ with respect to $\rho$ and $\theta$  respectively, while by exploiting~\eqref{e:pol} we get that $e_{\theta} =  R / (\gamma - 1)$. 
Also, to simplify the notations we write 
\begin{equation}
\label{e:hyp:a22}
    A_{21} = \frac{1}{ \theta}
    \left(
    \begin{array}{ccc}
               p_{\rho}  \\
               0  \\
    \end{array}
    \right)
    \qquad 
    A_{22} = \frac{1}{ \theta }
     \left(
    \begin{array}{ccc}
              \rho v -  \nu' \rho_x &  p_{\theta} \\
                      p_{\theta}  - 
                       \nu v_x /  \theta  
                      &    \rho v e_{\theta} / \theta  
                       -  k' \rho_x  /  \theta  \\              \end{array}
    \right)  
\end{equation}
and 
$$
     b = \frac{1}{ \theta }
    \left(
    \begin{array}{ccc}
            \nu   & 0  \\
            0 &   k / \theta  \\ 
    \end{array}
    \right) 
    \qquad 
    a_{11} = p_{\rho}  / (\theta \rho ).  
$$
The condition of block linear degeneracy is violated here. To see this, let us focus on the case $\sigma =0$: the dimension of the kernel of $A_{11} = a_{11}  v $ is $0$ if $v \neq 0$, but it is $1$ when $v = 0$ (we recall that $p_\rho >0$). 

To see what in principle can go wrong, we focus on steady solutions, the situation for travelling waves being analogous. 
We set  
\begin{equation}
\label{e:h:ns:ux}
   w = \rho _x  \qquad \vec z = \Big( v_x, \, \theta_x \Big)^T.
\end{equation}
Then \eqref{e:steady:vispro} becames 
$$
      \left(
             \begin{array}{cc}
                   a_{11} v  & A_{21}^T(u) \\
                   A_{21}(u) & A_{22}(u, \, u_x) \\
             \end{array}
             \right) 
             \left(
             \begin{array}{ll}
                         w \\
                         \vec z \\
             \end{array}
             \right)
             = 
             \left(
             \begin{array}{cc}
                   0 & 0 \\
                   0 & b(u) \\
             \end{array}
             \right) 
             \left(
             \begin{array}{ll}
                         w_x \\
                         \vec z_x \\
             \end{array}
             \right), 	     
$$
i.e.
$$
  \left\{
\begin{array}{ll}
           a_{11}  v  w                 +    A_{21}^T \vec  z = 0 \\
            A_{21} w        +    A_{22} \vec z = b \vec z_x 
\end{array}           
\right.
$$
Assume $v \neq 0$, then  \eqref{e:steady:vispro} can be written as 
\begin{equation}
\label{e:h:ns:zx}
\left\{
\begin{array}{ll}
            w =  \displaystyle{    -  \frac{ A_{21}^T \vec{z} }{a_{11}   v}     } \\
            \vec z_x =   b^{-1} \displaystyle{\Big[     A_{22}   -   \frac{A_{21} A_{21}^T   }{ a_{11} v}   \Big] \vec z } \\
\end{array}
\right.
\end{equation}
Note that the matrix $b$ is invertible and hence the previous expression is well defined. 

By combining \eqref{e:h:ns:ux} and \eqref{e:h:ns:zx}, one obtains that the steady solutions of the Navier Stokes written in Eulerian coordinates satisfy a singular ordinary differential equation in the form 
\begin{equation}
\label{e:sin}
    \frac{d U}{ d x } = \frac{1}{ v} F(U)
\end{equation}
provided that $U = \Big( \rho, \, v, \, \theta, \, \vec z \Big)^T$ and 
$$
    F(U) =
    \left(
    \begin{array}{ccc}
                A_{21}^T \vec{z} / a_{11} \\
                v \, \vec z   \\
                 b^{-1} \displaystyle{\Big[     A_{22} v   -   A_{21} A_{21}^T  /a_{11}   \Big] \vec z } \\
    \end{array}
    \right)
$$
We say that equation~\eqref{e:sin} is \emph{singular} since in general $v$, the velocity of the fluid, can attain the value $0$.  

Let us now go back to the example in~\cite{BiaSpi:rie} we mentioned before: the example deals with a system with non continuously differentiable steady solutions. It turns out that, as a consequence of the fact that the system violates the block linear degeneracy, the steady solutions $V$ satisfy a singular ODE in the form
\begin{equation}
\label{e:odesi}
         \frac{dV}{dt} =  \frac{1}{  \zeta (V)}  G(V),
\end{equation}
where $\zeta$ is a real-valued, smooth function that can attain the value $0$. The loss of regularity experienced by $V$ is actually due to the fact that~\eqref{e:odesi} is singular. 

Summing up, we have the following: the condition of block linear degeneracy is satisfied by the compressible Navier Stokes equation written in Lagrangian coordinates, but it is violated by the same equation written in Eulerian coordinates. 
As a consequence, the viscous profiles of the equation written in Eulerian coordinates 
satisfy a singular ODE and hence might in principle experience a loss of regularity.  This suggests that we should look for more general conditions than the block linear degeneracy. Namely, we want to find a condition on the singular equation \eqref{e:odesi}
which is sufficiently weak to be satisfied by the Navier Stokes equation written in Eulerian coordinates and by other systems violating the block linear degeneracy. On the other side, this condition should be sufficiently strong to rule out any loss of regularity. 

The definition of this condition and the corresponding analysis of the singular equation~\eqref{e:odesi} is done in \cite{BiaSpi:ode}. The conditions defined there apply to a class of systems that do not satisfy the block linear degeneracy and, in particular, to the Navier Stokes in Eulerian coordinates. 

\begin{remark}
For a different approach to the analysis of the viscous profiles of the Navier Stokes equation in Eulerian coordinates we refer for example to Wagner \cite{Wagner} and to the references therein. 
\end{remark}
\bibliography{biblio}
\end{document}